\documentclass[11pt,twoside]{amsart}
\usepackage{graphicx}
\usepackage{color}
\usepackage{amssymb}
\usepackage{amsmath}
\usepackage{amsthm}
\usepackage{amsfonts}
\usepackage{amsmath, amsfonts, amssymb}
\newtheorem{theorem}{Theorem}[section]

\newtheorem{de}{Definition}
\newtheorem{lm}{Lemma}

\newtheorem{theor}{Theorem}

\newtheorem{corollary}[theorem]{Corollary}

\numberwithin{equation}{section}

\begin{document}
			\title{A Complete classification of two-dimensional endo-commutative algebras over an arbitrary field}
		\author{D.Asrorov, U.Bekbaev, I.Rakhimov}
		
		\thanks{{\scriptsize
				emails: $^1$96asrorovdiorjon@mail.ru; $^2$uralbekbaev@gmail.com; $^3$isamiddin@uitm.edu.my}}
		\maketitle
		\begin{center}
			\address{$^1$Samarqand State University, Samarqand, Uzbekistan\\ $^2$Turin Polytechnic University in Tashkent, Tashkent, Uzbekistan\\ $^3$Universiti Teknologi MARA (UiTM), Shah Alam, Malaysia }
			\end{center}
		
			\begin{abstract}
		In the paper we describe the class of all 2-dimensional endo-commutative algebras over any basic field. Thereby, we extend S.-E. Takahasi, K. Shirayanagi and M. Tsukada's recent
description of the class of all 2-dimensional endo-commutative algebras over $\mathbb{Z}_2$ to any basic field, where the concept of endo-commutative algebra was first introduced and motivations to study this class of algebras also had been given. We give canonical representatives of the isomorphism classes such algebras in dimension two over any basic field.
This elementary and self-contained exposition uses recent result one of the authors on
description of the category of all $2$-dimensional algebras over any basic field.
\end{abstract}


	
	\section{Introduction}
The classification problem, up to isomorphism, of a given class of algebras is one of the important and difficult problems of algebra. So far two approaches are known to the solution of the problem. One of them is structural (basis free, invariant) approach. For instance, the classification of finite dimensional simple and semi-simple associative algebras by Wedderburn and simple and semi-simple Lie algebras by Cartan are examples of such approach. But it is observed that this approach becomes more difficult when one considers more general types of algebras. Another approach to the solution of the problem is coordinate based. Many researchers have used this type of approach to classify various, mainly finite dimensional, classes of algebras: associative \cite{KSTT,Mazzola}, Lie \cite{Jac, Mor, Mur1, Mur2, Mur3}, Jordan \cite{Kash} and Leibniz \cite{AOR, Casas, Demir, KhudRK}. These two approaches are somehow complementary to each other.

There were attempts to classify all fixed-dimensional algebras, for example, in \cite{P2000} a classification of all $2$-dimensional algebras, by basis free approach, was stated over any basic field. One disadvantage of the basis free approach is that an application of the obtained classification result to classification of a given class of algebras is hardly possible. In this respect the coordinate based classification has advantage over it. For the coordinated based approach in classification all $2$-dimensional algebras, over fields with some constraints, one can see \cite{ABR,G2011,IK} and in \cite{LJM2,SI} its different applications.


The concept of endo-commutative algebra was first introduced in \cite{TST}, where the authors gave a complete classification of two-dimensional endo-commutative algebras over the field of two elements. As well as, the paper contains a justification to study the class of endo-commutative algebras. In the present paper we give a complete classification of all endo-commutative algebras structures on two-dimensional vector space over any basic. The result of the paper based on a result obtained recently, in \cite{UralBekbaev} on complete classification of all $2$-dimensional algebras over any basic field. The result of the paper \cite{TST} comes here as a particular case.

The organization of the paper is as follows. The next section is Preliminaries, where we include the necessary definitions and results to be used throughout the paper. The main results of the paper are in Section 3 and onward. Sections 3.1, 3.2 and 3.3 contain the description of all two-dimensional endo-commutative algebras over a field $\mathbb{F}$ of the characteristic is not 2,3, the characteristic 2 and the characteristic 3, respectively. Here, in Corollary 3.1 we recover the result obtained in \cite{TST}. Sections 4 contains the classification of all two-dimensional curled algebras over any basic field and Section 5 is devoted to the description of two-dimensional endo-commutative curled algebras, where we show how to recover the result of \cite{TST} on the endo-commutative curled algebras.

\section{Preliminaries}
We begin by shortly recalling some concepts which will be using in the paper.

Let $\mathbf{A}$ be an $n$-dimensional algebra over a field $\mathbb{F}$ and $\mathbf{e}=\{\mathrm{e}_1,\mathrm{e}_2,...,\mathrm{e}_n\}$ be a basis of the underlying vector space of $\mathbf{A}$. Then on the basis $\mathbf{e}$ the algebra $\mathbf{A}$ is represented by a $n \times n^2$ matrix (called the matrix of structure constant, shortly MSC) $$A=\left(\begin{array}{ccccccccccccc}a_{11}^1&a_{12}^1&...&a_{1n}^1&a_{21}^1&a_{22}^1&...&a_{2n}^1&...&a_{n1}^1&a_{n2}^1&...&a_{nn}^1\\ a_{11}^2&a_{12}^2&...&a_{1n}^2&a_{21}^2&a_{22}^2&...&a_{2n}^2&...&a_{n1}^2&a_{n2}^2&...&a_{nn}^2 \\
...&...&...&...&...&...&...&...&...&...&...&...&...\\ a_{11}^n&a_{12}^n&...&a_{1n}^n&a_{21}^n&a_{22}^n&...&a_{2n}^n&...&a_{n1}^n&a_{n2}^n&...&a_{nn}^n\end{array}\right)$$ as follows
$$\mathrm{e}_i \cdot \mathrm{e}_j=\sum\limits_{k=1}^n a_{ij}^k\mathrm{e}_k, \ \mbox{where}\ i,j=1,2,...,n.$$
Therefore, the product on $\mathbf{A}$ with respect to the basis $\mathbf{e}$ is written as follows
\begin{equation} \label{Product}
	\mathrm{x}\cdot \mathrm{y}=\mathbf{e} A(x\otimes y)
\end{equation}
 for any $\mathrm{x}=\mathbf{e}x,\mathrm{y}=\mathbf{e}y,$
	where $x=(x_1, x_2,...,x_n)^T,$ and  $y=(y_1, y_2,...,y_n)^T$ are column coordinate vectors of $\mathrm{x}$ and $\mathrm{y},$ respectively, $x\otimes y$ is the tensor(Kronecker)
 product of the vectors $x$ and $y$. Now and onward for the product ``$\mathrm{x}\cdot \mathrm{y}$'' on $\mathbf{A}$ we use the juxtaposition ``$\mathrm{x} \mathrm{y}$''.
 	
 	If $\mathbf{A}$ is a two-dimensional algebra over a field $\mathbb{F}$ and $\mathbf{e}=\{\mathrm{e}_1, \mathrm{e}_2\}$ is a basis then MSC of $\mathbf{A}$ looks like	
 	 \[A=\left(\begin{array}{cccc} \alpha_1 & \alpha_2 & \alpha_3 &\alpha_4\\ \beta_1 & \beta_2 & \beta_3 &\beta_4\end{array}\right).\] In \cite{UralBekbaev} the following classification, up to isomorphism, of two-dimensional algebras, over any basic field $\mathbb{F}$, in terms of their MSC, was given.
 	
 	 	\begin{theor} \label{Char0} Any non-trivial $2$-dimensional algebra over a field $\mathbb{F}$ $(Char(\mathbb{F})\neq 2,3)$ is isomorphic to only one of the following listed, by their matrices of structure constants, such algebras:
 		\begin{itemize}
 			\item $A_1(\mathrm{c})=\begin{pmatrix}
 			\alpha_1&\alpha_2&1+\alpha_2&\alpha_4 \\
 			\beta_1& -\alpha_1& 1-\alpha_1& -\alpha_2
 			\end{pmatrix},$
 			where $ \mathrm{c}=(\alpha_1,\alpha_2,\alpha_4,\beta_1)\in\mathbb{F}^4$
 			
 			\item $A_2(\mathrm{c})=\begin{pmatrix}
 			\alpha_1&0&0&\alpha_4 \\
 			1& \beta_2& 1-\alpha_1& 0
 			\end{pmatrix},$
 			where $ \mathrm{c}=(\alpha_1,\alpha_4,\beta_2)\in\mathbb{F}^3$ and $\alpha_4\neq0$
 			
 			\item $A_3(\mathrm{c})=\begin{pmatrix}
 			\alpha_1&0&0&\alpha_4 \\
 			0& \beta_2& 1-\alpha_1& 0
 			\end{pmatrix}\simeq\begin{pmatrix}
 			\alpha_1&0&0&a^2\alpha_4 \\
 			0& \beta_2& 1-\alpha_1& 0
 			\end{pmatrix},$  where $ \mathrm{c}=(\alpha_1,\alpha_4,\beta_2)\in\mathbb{F}^3,$ $a\in\mathbb{F}$ and $ a\neq 0$
 			
 			\item $A_4(\mathrm{c})=\begin{pmatrix}
 			0&1&1&0 \\
 			\beta_1& \beta_2& 1& -1
 			\end{pmatrix},$  where $ \mathrm{c}=(\beta_1,\beta_2)\in\mathbb{F}^2$
 			
 			\item $A_5(\mathrm{c})=\begin{pmatrix}
 			\alpha_1&0&0&0 \\
 			1& 2\alpha_1-1& 1-\alpha_1& 0
 			\end{pmatrix},$\ \ where $ \mathrm{c}=\alpha_1\in\mathbb{F}$
 			
 			\item $A_6(\mathrm{c})=\begin{pmatrix}
 			\alpha_1&0&0&\alpha_4 \\
 			1& 1-\alpha_1& -\alpha_1& 0
 			\end{pmatrix},$ where $ \mathrm{c}=(\alpha_1,\alpha_4)\in\mathbb{F}^2$ and $ \alpha_4\neq0$
 			
 			\item $A_7(\mathrm{c})=\begin{pmatrix}
 			\alpha_1&0&0&\alpha_4 \\
 			0& 1-\alpha_1& -\alpha_1& 0
 			\end{pmatrix}\simeq\begin{pmatrix}
 			\alpha_1&0&0&a^2\alpha_4 \\
 			0& 1-\alpha_1& -\alpha_1& 0
 			\end{pmatrix},$ where $ \mathrm{c}=(\alpha_1,\alpha_4)\in\mathbb{F}^2, $ $ a\in\mathbb{F}$ and $ a\neq 0$
 			
 			\item $A_8(\mathrm{c})=\begin{pmatrix}
 			0&1&1&0 \\
 			\beta_1& 1&0& -1
 			\end{pmatrix},$ where $ \mathrm{c}=\beta_1\in\mathbb{F}$
 			
 			\item $A_9=\begin{pmatrix}
 			\frac{1}{3}&0&0&0\\
 			1&\frac{2}{3}&-\frac{1}{3}&0
 			\end{pmatrix},
 			$
 			
 			\item $A_{10}(\mathrm{c})=\begin{pmatrix}
 			0&1&1&1\\
 			\beta_1&0&0&-1
 			\end{pmatrix}\simeq\begin{pmatrix}
 			0&1&1&1\\
 			\beta_1^{'}(a)&0&0&-1
 			\end{pmatrix},$\\ where $ \mathrm{c}=\beta_1\in\mathbb{F},$ the polynomial $(\beta_1t^3-3t-1)(\beta_1t^2+\beta_1t+1)(\beta_1^2t^3+6\beta_1t^2+3\beta_1t+\beta_1-2)$ has no root in $\mathbb{F}$,
 			$a\in\mathbb{F}$ and $ \beta_1^{'}(t)=\frac{(\beta_1^2t^3+6\beta_1t^2+3\beta_1t+\beta_1-2)^2}{(\beta_1t^2+\beta_1t+1)^3}$
 			
 			\item $A_{11}(\mathrm{c})=\begin{pmatrix}
 			0&0&0&1 \\
 			\beta_1&0&0&0\\
 			\end{pmatrix}\simeq\begin{pmatrix}
 			0&0&0&1 \\
 			a^{3}\beta_1^{\pm1}&0&0&0\\
 			\end{pmatrix},$ where the polynomial $ \beta_1-t^3$ has no root in $\mathbb{F},$\ $ a, \mathrm{c}=\beta_1\in \mathbb{F}$  and $ a, \beta_1\neq0$
 			
 			\item $A_{12}(\mathrm{c})=\begin{pmatrix}
 			0&1&1&0 \\
 			\beta_1&0&0&-1
 			\end{pmatrix}\simeq\begin{pmatrix}
 			0&1&1&0 \\
 			a^2\beta_1&0&0&-1
 			\end{pmatrix},$
 			where $a, \mathrm{c}=\beta_1\in\mathbb{F}$ and $ a\neq 0$
 			
 			\item $A_{13}=\begin{pmatrix}
 			0&0&0&0 \\
 			1&0&0&0
 			\end{pmatrix}.
 			$
 		\end{itemize}
 	\end{theor}
 \begin{theor} \label{Char2} Any non-trivial $2$-dimensional algebra over a field $\mathbb{F}$ $(Char(\mathbb{F})=2)$ is isomorphic to only one of the following listed by their matrices of structure constants, such algebras:
 	\begin{itemize}
 		\item
 		$ A_{1,2}(\mathrm{c})=\begin{pmatrix}
 		\alpha_1&\alpha_2&1+\alpha_2&\alpha_4 \\
 		\beta_1& \alpha_1& 1+\alpha_1&\alpha_2
 		\end{pmatrix},
 		$ where $ \mathrm{c}=(\alpha_1,\alpha_2,\alpha_4,\beta_1)\in\mathbb{F}^4$
 		
 		\item
 		$ A_{2,2}(\mathrm{c})=\begin{pmatrix}
 		\alpha_1&0&0&\alpha_4 \\
 		1& \beta_2& 1+\alpha_1& 0
 		\end{pmatrix},
 		$ where $ \mathrm{c}=(\alpha_1,\alpha_4,\beta_2)\in\mathbb{F}^3$ and $\alpha_4\neq 0$
 		\item
 		$ A_{2,2}(\alpha_1,0,1)=\begin{pmatrix}
 		\alpha_1&0&0&0 \\
 		1& 1& 1+\alpha_1& 0
 		\end{pmatrix},
 		$ where $\alpha_1\in\mathbb{F}$
 		
 		\item $ A_{3,2}(\mathrm{c})=\begin{pmatrix}
 		\alpha_1&0&0&\alpha_4 \\
 		0& \beta_2& 1+\alpha_1& 0
 		\end{pmatrix}\simeq\begin{pmatrix}
 		\alpha_1&0&0&a^2\alpha_4 \\
 		0& \beta_2& 1+\alpha_1& 0
 		\end{pmatrix}
 		,$ where $ \mathrm{c}=(\alpha_1,\alpha_4,\beta_2)\in\mathbb{F}^3,$ $a\in\mathbb{F}$ and $ a\neq 0$

 		\item $ A_{4,2}(\mathrm{c})=\begin{pmatrix}
 		\alpha_1&1&1&0 \\
 		\beta_1& \beta_2& 1+\alpha_1& 1
 		\end{pmatrix}\simeq\begin{pmatrix}
 		\alpha_1&1&1&0 \\
 		\beta_1+(1+\beta_2)a+a^2& \beta_2& 1+\alpha_1& 1
 		\end{pmatrix},
 		$\\ where $ \mathrm{c}=(\alpha_1,\beta_1,\beta_2)\in\mathbb{F}^3$
 		
 		\item $A_{5,2}(\mathrm{c})=\begin{pmatrix}
 		\alpha_1&0&0&\alpha_4 \\
 		1& 1+\alpha_1& \alpha_1& 0
 		\end{pmatrix},
 		$ where $ \mathrm{c}=(\alpha_1,\alpha_4)\in\mathbb{F}^2$ and $\alpha_4\neq0$
 		\item $A_{5,2}(1,0)=\begin{pmatrix}
 		1&0&0&0 \\
 		1&0&1& 0
 		\end{pmatrix},
 		$
 		
 		\item $A_{6,2}(\mathrm{c})=\begin{pmatrix}
 		\alpha_1&0&0&\alpha_4 \\
 		0& 1+\alpha_1& \alpha_1& 0
 		\end{pmatrix}\simeq\begin{pmatrix}
 		\alpha_1&0&0&a^2\alpha_4 \\
 		0& 1+\alpha_1& \alpha_1& 0
 		\end{pmatrix}
 		,$ where $ \mathrm{c}=(\alpha_1,\alpha_4)\in\mathbb{F}^2,$ $a\in \mathbb{F}$ and $a\neq0$
 		
 		\item $A_{7,2}(\mathrm{c})=\begin{pmatrix}
 		\alpha_1&1&1&0 \\
 		\beta_1& 1+\alpha_1& \alpha_1& 1
 		\end{pmatrix}\simeq\begin{pmatrix}
 		\alpha_1&1&1&0 \\
 		\beta_1+a\alpha_1+a+a^2& 1+\alpha_1& \alpha_1& 1
 		\end{pmatrix}
 		,$ where $ \mathrm{c}=(\alpha_1,\beta_1)\in\mathbb{F}^2$ and $a\in\mathbb{F}$
 		
 		\item $A_{8,2}(\mathrm{c})=\begin{pmatrix}
 		0&1&1&1 \\
 		\beta_1& 0&0& 1
 		\end{pmatrix}\simeq\begin{pmatrix}
 		0&1&1&1 \\
 		\beta_1^{'}(a)& 0&0& 1
 		\end{pmatrix}
 		,$ where the polynomial\\ $(\beta_1t^3+t+1)(\beta_1t^2+\beta_1t+1)$ has no root in $\mathbb{F},$
 		$a\in\mathbb{F}$ and $ \beta_1^{'}(t)=\frac{(\beta_1^2t^3+\beta_1t+\beta_1)^2}{(\beta_1t^2+\beta_1t+1)^3}$
 		
 		\item $A_{9,2}(\mathrm{c})=\begin{pmatrix}
 		0&0&0&1\\
 		\beta_1&0&0&0
 		\end{pmatrix}\simeq\begin{pmatrix}
 		0&0&0&1\\
 		a^3\beta_1^{\pm1}&0&0&0
 		\end{pmatrix}
 		,$ where $ a, \mathrm{c}=\beta_1\in\mathbb{F}$ and $a\neq 0,$
 		
 		the polynomial $ \beta_1+t^3$ has no root in $\mathbb{F}$
 		
 		\item $A_{10,2}(\mathrm{c})=\begin{pmatrix}
 		1&1&1&0\\
 		\beta_1&1&1&1
 		\end{pmatrix}\simeq\begin{pmatrix}
 		1&1&1&0\\
 		\beta_1+a+a^2&1&1&1
 		\end{pmatrix}
 		,$ where $a, \mathrm{c}=\beta_1\in\mathbb{F}$
 		
 		\item $A_{11,2}(\mathrm{c})=\begin{pmatrix}
 		0&1&1&0 \\
 		\beta_1&0&0&1\\
 		\end{pmatrix}\simeq\begin{pmatrix}
 		0&1&1&0 \\
 		b^2(\beta_1+a^2)&0&0&1\\
 		\end{pmatrix}
 		,$ where $a, b\in\mathbb{F}$ and $ b\neq 0$
 		
 		\item $A_{12,2}=\begin{pmatrix}
 		0&0&0&0 \\
 		1&0&0&0\\
 		\end{pmatrix}
 		$
  		
 	\end{itemize}
 \end{theor}

 \begin{theor} \label{Char3} Any non-trivial $2$-dimensional algebra over a field $\mathbb{F}$ $(Char(\mathbb{F})=3)$ is isomorphic to only one of the following, listed by their matrices of structure constants, such algebras:
 	\begin{itemize}
 		\item $A_{1,3}(\mathrm{c})=\begin{pmatrix}
 		\alpha_1&\alpha_2&\alpha_2+1&\alpha_4 \\
 		\beta_1& -\alpha_1& 1-\alpha_1& -\alpha_2
 		\end{pmatrix},$ where $ \mathrm{c}=(\alpha_1,\alpha_2,\alpha_4,\beta_1)\in\mathbb{F}^4$
 		
 		\item $A_{2,3}(\mathrm{c})=\begin{pmatrix}
 		\alpha_1&0&0&\alpha_4 \\
 		1& \beta_2& 1-\alpha_1& 0
 		\end{pmatrix},$ where $ \mathrm{c}=(\alpha_1,\alpha_4,\beta_2)\in\mathbb{F}^3,$ and $\alpha_4\neq0$
 		
 		\item $A_{3,3}(\mathrm{c})=\begin{pmatrix}
 		\alpha_1&0&0&\alpha_4 \\
 		0& \beta_2& 1-\alpha_1& 0
 		\end{pmatrix}\simeq\begin{pmatrix}
 		\alpha_1&0&0&a^2\alpha_4 \\
 		0& \beta_2& 1-\alpha_1& 0
 		\end{pmatrix},$ where $ \mathrm{c}=(\alpha_1,\alpha_4,\beta_2)\in\mathbb{F}^3,$ $ a\in\mathbb{F}$ and $ a\neq 0$
 		
 		\item $A_{4,3}(\mathrm{c})=\begin{pmatrix}
 		0&1&1&0 \\
 		\beta_1& \beta_2& 1& -1
 		\end{pmatrix},$ where $ \mathrm{c}=(\beta_1,\beta_2)\in\mathbb{F}^2$
 		
 		\item $A_{5,3}(\mathrm{c})=\begin{pmatrix}
 		\alpha_1&0&0&0 \\
 		1& 2\alpha_1-1& 1-\alpha_1& 0
 		\end{pmatrix},$ where $ \mathrm{c}=\alpha_1\in\mathbb{F}$
 		
 		\item $A_{6,3}(\mathrm{c})=\begin{pmatrix}
 		\alpha_1&0&0&\alpha_4 \\
 		1& 1-\alpha_1& -\alpha_1& 0
 		\end{pmatrix},$ where $ \mathrm{c}=(\alpha_1,\alpha_4)\in\mathbb{F}^2$ and $ \alpha_4\neq0$
 		
 		\item $A_{7,3}(\mathrm{c})=\begin{pmatrix}
 		\alpha_1&0&0&\alpha_4 \\
 		0& 1-\alpha_1& -\alpha_1& 0
 		\end{pmatrix}\simeq\begin{pmatrix}
 		\alpha_1&0&0&a^2\alpha_4 \\
 		0& 1-\alpha_1& -\alpha_1& 0
 		\end{pmatrix},$ where $ \mathrm{c}=(\alpha_1,\alpha_4)\in\mathbb{F}^2,$ $a\in\mathbb{F}$ and $ a\neq 0$
 		
 		\item $A_{8,3}(\mathrm{c})=\begin{pmatrix}
 		0&1&1&0 \\
 		\beta_1& 1&0& -1
 		\end{pmatrix},$ where $ \mathrm{c}=\beta_1\in\mathbb{F}$
 		
 		\item $A_{9,3}(\beta_1)=\begin{pmatrix}
 		0&1&1&1\\
 		\beta_1&0&0&-1
 		\end{pmatrix}\simeq\begin{pmatrix}
 		0&1&1&1\\
 		\beta_1'(a)&0&0&-1
 		\end{pmatrix},$ where the polynomial \\ $(\beta_1-t^3)(\beta_1t^2+\beta_1t+1)(\beta_1^2t^3+\beta_1-2)$ has no root in $\mathbb{F},$
 		$a\in\mathbb{F}$ and $ \beta_1^{'}(t)=\frac{(\beta_1^2t^3+\beta_1-2)^2}{(\beta_1t^2+\beta_1t+1)^3}$
 		
 		\item $A_{10,3}(\mathrm{c})=\begin{pmatrix}
 		0&0&0&1\\
 		\beta_1&0&0&0
 		\end{pmatrix}\simeq\begin{pmatrix}
 		0&0&0&1\\
 		a^{3}\beta_1^{\pm1}&0&0&0
 		\end{pmatrix},$\\  where the polynomial $ \beta_1-t^3$ has no root, $ a, \mathrm{c}=\beta_1\in\mathbb{F}$ and  $ a, \beta_1\neq 0$
 		
 		\item $A_{11,3}(\mathrm{c})=\begin{pmatrix}
 		0&1&1&0 \\
 		\beta_1&0&0&-1\\
 		\end{pmatrix}\simeq\begin{pmatrix}
 		0&1&1&0 \\
 		a^{2}\beta_1&0&0&-1\\
 		\end{pmatrix},
 		$ where  $a, \mathrm{c}=\beta_1\in\mathbb{F},$ $a\neq 0$
 		
 		\item $A_{12,3}=\begin{pmatrix}
 		1&0&0&0 \\
 		1&-1&-1&0
 		\end{pmatrix}
 		,$
 		
 		\item $A_{13,3}=\begin{pmatrix}
 		0&0&0&0 \\
 		1&0&0&0
 		\end{pmatrix}.$
 	\end{itemize}
 \end{theor}

%

 \begin{de} An algebra $\mathbf{A}$ is said to be endo-commutative if $\mathrm{x}^2\mathrm{y}^2=(\mathrm{x} \mathrm{y})^2,$ for any $\mathrm{x}, \mathrm{y} \in \mathbf{A}.$
 		
 \end{de}
	\begin{lm} An algebra $\mathbf{A}$ is endo-commutative if and only if \begin{equation} \label{ECA} A(A\otimes A)(x^{\otimes 2}\otimes y^{\otimes 2}-(x\otimes y)^{\otimes 2})=0,\end{equation} where $A$ is MSC of $\mathbf{A}$.
 		
 \end{lm}
 \begin{proof}
	According to (\ref{Product}) we write $$\mathrm{x}^2=\mathbf{e}Ax^{\otimes 2}, \ \ \mathrm{y}^2=\mathbf{e}Ay^{\otimes 2}, \ \ \mathrm{x}^2\mathrm{y}^2=\mathbf{e}A(Ax^{\otimes 2}\otimes Ay^{\otimes 2})=\mathbf{e}A(A\otimes A)(x^{\otimes 2}\otimes y^{\otimes 2})$$ and
	$(\mathrm{x} \mathrm{y})^2=\mathbf{e}A(A(x\otimes y)\otimes A(x\otimes y))=\mathbf{e}A(A\otimes A)(x\otimes y)^{\otimes 2}.$
	Therefore, the equality $\mathrm{x}^2\mathrm{y}^2=(\mathrm{x} \mathrm{y})^2$ in terms of MSC and the coordinate vectors is written as follows

	$\mathbf{e}A(A\otimes A)(x^{\otimes 2}\otimes y^{\otimes 2})=\mathbf{e}A(A\otimes A)(x\otimes y)^{\otimes 2}$ i.e.,
	$A(A\otimes A)(x^{\otimes 2}\otimes y^{\otimes 2}-(x\otimes y)^{\otimes 2})=0$\\
what is required to get.
\end{proof}

In \cite{TST} the class of endo-commutative algebras was split into two: curled and straight algebras. There one can find a list of such curled algebras, up to isomorphism, over the field $ \mathbb{F}=\mathbb{Z}_2$ as well. The definition was given as follows.
	 \begin{de} An algebra $\mathbf{A}$ is said to be curled if $\mathrm{x}^2=\mathrm{\lambda(x)}\mathrm{x},$ for any $\mathrm{x}\in \mathbf{A},$ where $ \mathrm{\lambda(x)}\in \mathbb{F}.$
 \end{de}

Let $\mathbf{e}=\{\mathrm{e}_1, \mathrm{e}_2,..., \mathrm{e}_n\}$ be a basis of $\mathbf{A}$. According to (\ref{Product}) one has $\mathrm{x}^2=\mathbf{e}Ax^{\otimes 2},$ $\mathrm{\lambda(x)}\mathrm{x}=\mathbf{e}\lambda(x)x$ and therefore, the equality $\mathrm{x}^2=\mathrm{\lambda(x)}\mathrm{x}$ in terms of MSC of $A$ and the coordinate vector of $\mathrm{x}$ can be written as follows \begin{equation}
\label{C} Ax^{\otimes 2}-\lambda(x)x=0.\end{equation}

In Section  \ref{CCA} we give a complete classification of all two-dimensional curled algebras over any basic field and then split them into curled and straight algebras. The lists of these algebras given in \cite{TST} come as a particular case.

First we start with the classification of all two-dimensional endo-commutative algebras based on the result of Theorems \ref{Char0}, \ref{Char2} and \ref{Char3}.

Let now $\mathbf{A}$ be a two-dimensional algebra over a field $\mathbb{F}$ and	
let \[A=\left(\begin{array}{cccc} \alpha_1 & \alpha_2 & \alpha_3 &\alpha_4\\ \beta_1 & \beta_2 & \beta_3 &\beta_4\end{array}\right)\] be its MSC on a basis $\mathbf{e}=\{\mathrm{e}_1, \mathrm{e}_2\}$, $x=(x_1,x_2)$ and $y=(y_1,y_2)$. Then (\ref{ECA})	is nothing else
	\begin{equation} \label{IDAB} \left\lbrace\begin{array}{llll}
x_1^2y_1y_2(A_3-A_5)+x_1^2y_2^2(A_4-A_6)+x_1x_2y_1y_2(A_6-A_4+A_{11}-A_{13})\\
\hfill +x_1x_2y_1^2(A_5-A_3)+x_1x_2y_2^2(A_{12}-A_{14})+x_2^2y_1^2(A_{13}-A_{11})+x_2^2y_1y_2(A_{14}-A_{12})&=&0\\
x_1^2y_1y_2(B_3-B_5)+x_1^2y_2^2(B_4-B_6)+x_1x_2y_1y_2(B_6-B_4+B_{11}-B_{13})
\\
\hfill+x_1x_2y_1^2(B_5-B_3)+x_1x_2y_2^2(B_{12}-B_{14})+x_2^2y_1^2(B_{13}-B_{11})+x_2^2y_1y_2(B_{14}-B_{12})&=&0,
\end{array}\right.
\end{equation}
where $$\begin{array}{llllllllll} A_3&=& \alpha_1^2\alpha_3+\alpha_1\alpha_2\beta_3+\alpha_3^2\beta_1+\alpha_4\beta_1\beta_3 \ \
&A_4&=&\alpha_1^2\alpha_4+\alpha_1\alpha_2\beta_4+\alpha_3\alpha_4\beta_1+\alpha_4\beta_1\beta_4 \\
 A_5&=&\alpha_1^2\alpha_2+\alpha_2^2\beta_1+\alpha_1\alpha_3\beta_2+\alpha_4\beta_1\beta_2 \ \
 &A_6&=&\alpha_1\alpha_2^2+\alpha_2^2\beta_2+\alpha_2\alpha_3\beta_2+\alpha_4\beta_2^2 \\
 A_{11}&=&\alpha_1\alpha_3^2+\alpha_2\alpha_3\beta_3+\alpha_3^2\beta_3+\alpha_4\beta_3^2 \ \
& A_{12}&=&\alpha_1\alpha_3\alpha_4+\alpha_2\alpha_3\beta_4+\alpha_3\alpha_4\beta_3+\alpha_4\beta_3\beta_4 \\
 A_{13}&=&\alpha_1^2\alpha_4+\alpha_2\alpha_4\beta_1+\alpha_1\alpha_3\beta_4+\alpha_4\beta_1\beta_4 \ \
& A_{14}&=&\alpha_1\alpha_2\alpha_4+\alpha_2\alpha_4\beta_2+\alpha_2\alpha_3\beta_4+\alpha_4\beta_2\beta_4 \\
 B_3&=& \alpha_1\alpha_3\beta_1+\alpha_1\beta_2\beta_3+\alpha_3\beta_1\beta_3+\beta_1\beta_3\beta_4 \ \
 &B_4&=&\alpha_1\alpha_4\beta_1+\alpha_1\beta_2\beta_4+\alpha_4\beta_1\beta_3+\beta_1\beta_4^2 \\
 B_5&=&\alpha_1\alpha_2\beta_1+\alpha_2\beta_1\beta_2+\alpha_1\beta_2\beta_3+\beta_1\beta_2\beta_4 \ \
 &B_6&=&\alpha_2^2\beta_1+\alpha_2\beta_2^2+\alpha_2\beta_2\beta_3+\beta_2^2\beta_4 \\
 B_{11}&=&\alpha_3^2\beta_1+\alpha_3\beta_2\beta_3+\alpha_3\beta_3^2+\beta_3^2\beta_4 \ \
 &B_{12}&=&\alpha_3\alpha_4\beta_1+\alpha_3\beta_2\beta_4+\alpha_4\beta_3^2+\beta_3\beta_4^2 \\
 B_{13}&=&\alpha_1\alpha_4\beta_1+\alpha_4\beta_1\beta_2+\alpha_1\beta_3\beta_4+\beta_1\beta_4^2 \ \
& B_{14}&=&\alpha_2\alpha_4\beta_1+\alpha_4\beta_2^2+\alpha_2\beta_3\beta_4+\beta_2\beta_4^2.
\end{array}$$

Note that the set of functions $\{x_1^2y_1y_2,\ x_1^2y_2^2,\ x_1x_2y_1^2,\ x_1x_2y_1y_2,\ x_1x_2y_2^2,\ x_2^2y_1^2,\ x_2^2y_1y_2\}$ is linearly independent. Therefore, the system (\ref{IDAB}) in terms of $A_i$ and $B_j$ ($i,j=3,4,...,13$) can be rewritten as follows
\begin{equation} \label{IDABK} \left\lbrace\begin{array}{rrrrrrrrrrr}
A_3-A_5&=&0\ \ \
A_4-A_6&=&0\ \ \
A_{12}-A_{14}&=&0\ \ \
A_{13}-A_{11}&=&0\\
B_3-B_5&=&0\ \ \
B_4-B_6&=&0\ \ \
B_{12}-B_{14}&=&0\ \ \
B_{13}-B_{11}&=&0.
\end{array}\right.
\end{equation}

In terms of the structure constants $\alpha_i, \beta_j$ ($i,j=1,2,...,4$) the system (\ref{IDABK}) is written as follows

\begin{equation} \label{GE}
\left\lbrace\begin{array}{lccc}
(\alpha_1^2\alpha_3+\alpha_1\alpha_2\beta_3+\alpha_3^2\beta_1+\alpha_4\beta_1\beta_3)-(\alpha_1^2\alpha_2+\alpha_2^2\beta_1+\alpha_1\alpha_3\beta_2+\alpha_4\beta_1\beta_2)&=&0\\
(\alpha_1^2\alpha_4+\alpha_1\alpha_2\beta_4+\alpha_3\alpha_4\beta_1+\alpha_4\beta_1\beta_4)-(\alpha_1\alpha_2^2+\alpha_2^2\beta_2+\alpha_2\alpha_3\beta_2+\alpha_4\beta_2^2)&=&0\\
(\alpha_1\alpha_3\alpha_4+\alpha_2\alpha_3\beta_4+\alpha_3\alpha_4\beta_3+\alpha_4\beta_3\beta_4)-(\alpha_1\alpha_2\alpha_4+\alpha_2\alpha_4\beta_2+\alpha_2\alpha_3\beta_4+\alpha_4\beta_2\beta_4)&=&0\\
(\alpha_1^2\alpha_4+\alpha_2\alpha_4\beta_1+\alpha_1\alpha_3\beta_4+\alpha_4\beta_1\beta_4)-(\alpha_1\alpha_3^2+\alpha_2\alpha_3\beta_3+\alpha_3^2\beta_3+\alpha_4\beta_3^2)&=&0\\
(\alpha_1\alpha_3\beta_1+\alpha_1\beta_2\beta_3+\alpha_3\beta_1\beta_3+\beta_1\beta_3\beta_4)-(\alpha_1\alpha_2\beta_1+\alpha_2\beta_1\beta_2+\alpha_1\beta_2\beta_3+\beta_1\beta_2\beta_4)&=&0\\
(\alpha_1\alpha_4\beta_1+\alpha_1\beta_2\beta_4+\alpha_4\beta_1\beta_3+\beta_1\beta_4^2)-(\alpha_2^2\beta_1+\alpha_2\beta_2^2+\alpha_2\beta_2\beta_3+\beta_2^2\beta_4)&=&0\\
(\alpha_3\alpha_4\beta_1+\alpha_3\beta_2\beta_4+\alpha_4\beta_3^2+\beta_3\beta_4^2)-(\alpha_2\alpha_4\beta_1+\alpha_4\beta_2^2+\alpha_2\beta_3\beta_4+\beta_2\beta_4^2)&=&0\\
(\alpha_1\alpha_4\beta_1+\alpha_4\beta_1\beta_2+\alpha_1\beta_3\beta_4+\beta_1\beta_4^2)-(\alpha_3^2\beta_1+\alpha_3\beta_2\beta_3+\alpha_3\beta_3^2+\beta_3^2\beta_4)&=&0.
\end{array}\right.
\end{equation}
Further we refer to (\ref{GE}) as a general system of equations for two-dimensional endo-commutative algebras.

\section{A complete classification of $2$-dimensional endo-commutative algebras}
In this sections and the sections in row we give canonical representatives of the isomorphism classes of two-dimensional endo-commutative algebras over the of characteristic is not 2,3, the characteristic 2 and the characteristic 3, respectively.
\subsection{Classification over $\mathbb{F}$ with $\mathbf{Char(\mathbb{F})\neq2,3}$}\emph{}
The list of canonical representatives of the isomorphism classes of two-dimensional endo-commutative algebras over a field $\mathbb{F}$ of $(Char\mathbb{F}\neq 2,3)$ is given as follows.

\begin{theor} \label{char0Endo} Any non-trivial $2$-dimensional endo-commutative algebra over a field $\mathbb{F},$ $(Char(\mathbb{F})\neq 2,3)$ is isomorphic to only one of the following listed by their matrices of structure constants, such algebras:

\begin{itemize}
\item
$ A_3(\alpha_1,0,\beta_2)=\begin{pmatrix}
\alpha_1&0&0&0 \\
0& \beta_2& 1-\alpha_1& 0
\end{pmatrix}
,$ where $ \alpha_1,\beta_2\in\mathbb{F}$
\item
$ A_3(\frac{1}{2},\alpha_4,\frac{1}{2})=\begin{pmatrix}
\frac{1}{2}&0&0&\alpha_4 \\
0& \frac{1}{2}& \frac{1}{2}& 0
\end{pmatrix},
$ where $\alpha_4\in\mathbb{F}$ and $\alpha_4\neq0$
\item
$ A_3(\frac{1}{2},\alpha_4,-\frac{1}{2})=\begin{pmatrix}
\frac{1}{2}&0&0&\alpha_4 \\
0& -\frac{1}{2}& \frac{1}{2}& 0
\end{pmatrix},
$ where $\alpha_4\in\mathbb{F}$ and $\alpha_4\neq0$
\item $ A_5(\alpha_1)=\begin{pmatrix}
\alpha_1&0&0&0 \\
1& 2\alpha_1-1& 1-\alpha_1& 0
\end{pmatrix}
,$ where $ \alpha_1\in\mathbb{F}$
\item
$ A_7(\alpha_1,0)=\begin{pmatrix}
\alpha_1&0&0&0 \\
0& 1-\alpha_1& -\alpha_1& 0
\end{pmatrix}
,$ where $ \alpha_1\in\mathbb{F}$
\item $ A_7(\frac{1}{2},\alpha_4)=\begin{pmatrix}
\frac{1}{2}&0&0&\alpha_4 \\
0& \frac{1}{2}& -\frac{1}{2}& 0
\end{pmatrix},
$ where $\alpha_4\in\mathbb{F}$ and $\alpha_4\neq0$
\item $A_9=\begin{pmatrix}
\frac{1}{3}&0&0&0\\
1&\frac{2}{3}&-\frac{1}{3}&0
\end{pmatrix},
$

\item $ A_{10}(\beta_1)=\begin{pmatrix}
0&1&1&1\\
\beta_1&0&0&-1
\end{pmatrix}\simeq\begin{pmatrix}
0&1&1&1\\
\beta_1^{'}(a)&0&0&-1
\end{pmatrix},$\ \ where $ a, \beta_1\in\mathbb{F},$ the polynomial $(\beta_1t^3-3t-1)(\beta_1t^2+\beta_1t+1)(\beta_1^2t^3+6\beta_1t^2+3\beta_1t+\beta_1-2)$ has no root in $\mathbb{F},$
and $ \beta_1^{'}(t)=\frac{(\beta_1^2t^3+6\beta_1t^2+3\beta_1t+\beta_1-2)^2}{(\beta_1t^2+\beta_1t+1)^3}$

\item $A_{11}(\beta_1)=\begin{pmatrix}
	0&0&0&1 \\
	\beta_1^{\pm1}&0&0&0\\
\end{pmatrix}\simeq\begin{pmatrix}
	0&0&0&1 \\
	a^{3}\beta_1^{\pm1}&0&0&0\\
\end{pmatrix},$\ \ where the polynomial $ \beta_1-t^3$ has no root in $\mathbb{F},$ $a, \beta_1 \in\mathbb{F}$ and  $a\neq 0$

\item $A_{12}(\beta_1)=\begin{pmatrix}
	0&1&1&0 \\
	\beta_1&0&0&-1
\end{pmatrix}\simeq\begin{pmatrix}
0&1&1&0 \\
a^2\beta_1&0&0&-1
\end{pmatrix},$\ \
where $ a, \beta_1\in\mathbb{F}$ and $ a\neq 0$

\item $ A_{13}=\begin{pmatrix}
0&0&0&0 \\
1&0&0&0
\end{pmatrix}.
$
\end{itemize}
\end{theor}

\begin{proof}
To classify the two-dimensional endo-commutative algebras it suffices to solve the general system of equations (\ref{GE}) with respect to MSC of each $A_i,$ ($ i=1,2,...,13$) given in Theorem \ref{Char0}.

For
 $ A=A_1(\mathrm{c})=\begin{pmatrix}
\alpha_1&\alpha_2&1+\alpha_2&\alpha_4 \\
\beta_1& -\alpha_1& 1-\alpha_1& -\alpha_2
\end{pmatrix},
$
where $ \mathrm{c}=(\alpha_1,\alpha_2,\alpha_4,\beta_1)\in\mathbb{F}^4,$
the system of equations (\ref{GE}) looks like

$$ \left\lbrace\begin{array}{lllll}
(\alpha_1^2(1+\alpha_2)+\alpha_1\alpha_2(1-\alpha_1)+(1+\alpha_2)^2\beta_1+\alpha_4\beta_1(1-\alpha_1))\\ \hfill-(\alpha_1^2\alpha_2+\alpha_2^2\beta_1-\alpha_1^2(1+\alpha_2)-\alpha_1\alpha_4\beta_1)&=&0\\
(\alpha_1^2\alpha_4-\alpha_1\alpha_2^2+(1+\alpha_2)\alpha_4\beta_1-\alpha_2\alpha_4\beta_1)\\ \hfill -(\alpha_1\alpha_2^2-\alpha_1\alpha_2^2-\alpha_1\alpha_2(1+\alpha_2)+\alpha_4\alpha_1^2)&=&0\\
(\alpha_1(1+\alpha_2)\alpha_4+(1+\alpha_2)\alpha_4(1-\alpha_1)-\alpha_2\alpha_4(1-\alpha_1))\\ \hfill -(\alpha_1\alpha_2\alpha_4-\alpha_1\alpha_2\alpha_4+\alpha_4\alpha_1\alpha_2))&=&0\\
(\alpha_1^2\alpha_4+\alpha_2\alpha_4\beta_1+\alpha_1(1+\alpha_2)(-\alpha_2)-\alpha_4\beta_1\alpha_2)\\ \hfill-(\alpha_1(1+\alpha_2)^2+\alpha_2(1+\alpha_2)(1-\alpha_1)+(1+\alpha_2)^2(1-\alpha_1)+\alpha_4(1-\alpha_1)^2)&=&0\\
(\alpha_1(1+\alpha_2)\beta_1+(1+\alpha_2)\beta_1(1-\alpha_1)+\beta_1(1-\alpha_1)(-\alpha_2))\\ \hfill-(\alpha_1\alpha_2\beta_1-\alpha_1\alpha_2\beta_1+\beta_1(-\alpha_1)(-\alpha_2))&=&0\\
(\alpha_1\alpha_4\beta_1-\alpha_1\alpha_2\beta_2+\alpha_4\beta_1(1-\alpha_1)+\beta_1(-\alpha_2)^2)\\ \hfill-(\alpha_2^2\beta_1+\alpha_2(-\alpha_1)^2+\alpha_2(-\alpha_1)(1-\alpha_1)+(-\alpha_1)^2(-\alpha_2)&=&0\\
((1+\alpha_2)\alpha_4\beta_1+(1+\alpha_2)(-\alpha_1)(-\alpha_2)+\alpha_4(1-\alpha_1)^2+(1-\alpha_1)(-\alpha_2)^2)\\ \hfill-(\alpha_2\alpha_4\beta_1+\alpha_4(-\alpha_1)^2+\alpha_2(1-\alpha_1)(-\alpha_2)-\alpha_1(-\alpha_2)^2)&=&0\\
(\alpha_1\alpha_4\beta_1-\alpha_1\alpha_4\beta_1+\alpha_1(1-\alpha_1)(-\alpha_2)+\beta_1(-\alpha_2)^2)\\ \hfill-((1+\alpha_2)^2\beta_1+(1+\alpha_2)(-\alpha_1)(1-\alpha_1)+(1+\alpha_2)(1-\alpha_1)^2+(1-\alpha_1)^2(-\alpha_2))&=&0.\\
\end{array}\right.
$$

Simplifying we get an inconsistent system of equations. Thus,
%
%
%
%
%
%
%
%
%
%
%
%
%
%
%
%
%
%
%
%
%
%
there is no endo-commutative algebra among $A_1(\mathrm{c})$.

Similarly, if
$A= A_2(\mathrm{c})=\begin{pmatrix}
\alpha_1&0&0&\alpha_4 \\
1& \beta_2& 1-\alpha_1& 0
\end{pmatrix},
$ where $ \mathrm{c}=(\alpha_1,\alpha_4,\beta_2)\in\mathbb{F}^3$ and $\alpha_4\neq 0$,
then the system (\ref{GE}) implies $\alpha_4= 0$ which contradicts with $\alpha_4\neq 0$. Hence, $A_2(\mathrm{c})$ also does not contain an endo-commutative algebra.

%
%
%
%
%
%
%

Let consider
$A= A_3(\mathrm{c})=\begin{pmatrix}
\alpha_1&0&0&\alpha_4 \\
0& \beta_2& 1-\alpha_1& 0
\end{pmatrix},
$ where $ \mathrm{c}=(\alpha_1,\alpha_4,\beta_2)\in\mathbb{F}^3.$
Then the system of equations (\ref{GE}) for $A_3(\mathrm{c})$ looks like

$$ \left\lbrace\begin{array}{rrr}

\alpha_1^2\alpha_4-\alpha_4\beta_2^2&=&0\\

\alpha_1^2\alpha_4-(1-\alpha_1)^2\alpha_4&=&0\\

(1-\alpha_1)^2\alpha_4-\alpha_4\beta_2^2&=&0\\

\end{array}\right.
$$ which has the following solutions
: $(\alpha_1,0,\beta_2)$, $(\frac{1}{2},\alpha_4,-\frac{1}{2})$ and $(\frac{1}{2},\alpha_4,\frac{1}{2})$. These solutions produce the endo-commutative algebras:
$$ A_3(\alpha_1,0,\beta_2)=\begin{pmatrix}
\alpha_1&0&0&0 \\
0& \beta_2& 1-\alpha_1& 0
\end{pmatrix},
\ \  A_3\left(\frac{1}{2},\alpha_4,\frac{1}{2}\right)=\begin{pmatrix}
\frac{1}{2}&0&0&\alpha_4 \\
0& \frac{1}{2}& \frac{1}{2}& 0
\end{pmatrix}
$$
and
$$ A_3\left(\frac{1}{2},\alpha_4,-\frac{1}{2}\right)=\begin{pmatrix}
\frac{1}{2}&0&0&\alpha_4 \\
0& -\frac{1}{2}& \frac{1}{2}& 0
\end{pmatrix},
$$

where $ \alpha_1, \alpha_4, \beta_2\in\mathbb{F}$
and $\alpha_4\neq0$.
%
%

Letting $ A=A_4(\mathrm{c})=\begin{pmatrix}
0&1&1&0 \\
\beta_1& \beta_2& 1& -1
\end{pmatrix},
$ where $ \mathrm{c}=(\beta_1,\beta_2)\in\mathbb{F}^2,$ we get an inconsistent system of equations.
%
%
%
%
%
%

If
$ A=A_5(\mathrm{c})=\begin{pmatrix}
\alpha_1&0&0&0 \\
1& 2\alpha_1-1& 1-\alpha_1& 0
\end{pmatrix},
$ where $ \mathrm{c}=\alpha_1\in\mathbb{F},$ then all the equations of the system (\ref{GE}) become identities, therefore, all the algebras in this class are endo-commutative.

For
$ A=A_6(\mathrm{c})=\begin{pmatrix}
\alpha_1&0&0&\alpha_4 \\
1& 1-\alpha_1& -\alpha_1& 0
\end{pmatrix},
$ where $ \mathrm{c}=(\alpha_1,\alpha_4)\in\mathbb{F}^2$ and $\alpha_4\neq 0$, the system of equations (\ref{GE}) gives $\alpha_4=0$ which is a contradiction.
%
%
%
%
%

If
$A= A_7(\mathrm{c})=\begin{pmatrix}
\alpha_1&0&0&\alpha_4 \\
0& 1-\alpha_1& -\alpha_1& 0
\end{pmatrix},
 \ \mbox{where}\ \mathrm{c}=(\alpha_1,\alpha_4)\in\mathbb{F}^2$
then the system of equations (\ref{GE}) is equivalent to $\alpha_1^2\alpha_4-\alpha_4(1-\alpha_1)^2=0$ and we have the solutions:
$$(\alpha_1,0) \ \mbox{and}\ \left(\frac{1}{2},\alpha_4\right), \ \mbox{where}\ \alpha_4\neq 0,$$ i.e.,
$$A_7(\alpha_1,0)=\begin{pmatrix}
\alpha_1&0&0&0 \\
0& 1-\alpha_1& -\alpha_1& 0
\end{pmatrix}, \alpha_1\in \mathbb{F}\ \mbox{and} \ A_7\left(\frac{1}{2},\alpha_4\right)=\begin{pmatrix}
\frac{1}{2}&0&0&\alpha_4 \\
0& \frac{1}{2}& -\frac{1}{2}& 0
\end{pmatrix},\ \alpha_4\neq 0 \in \mathbb{F} $$ are endo-commutative.

Letting
\[ A=A_8(\mathrm{c})=\begin{pmatrix}
0&1&1&0 \\
\beta_1& 1&0& -1
\end{pmatrix},
\] where $ \mathrm{c}=\beta_1\in\mathbb{F}$
we obtain an inconsistent system of equations.

Verifying the system of equations (\ref{GE}) for $A_9,$ $A_{10}(\mathrm{c}),$ $A_{11}(\mathrm{c}),$ $A_{12}(\mathrm{c})$ and $A_{13}$ we get identities, i.e., all these algebras happen to be endo-commutative algebras.
\end{proof}
\subsection{Classification over $\mathbb{F}$ with $\mathbf{Char(\mathbb{F})= 2,3}$}\emph{}

The description of two-dimensional endo-commutative algebras over $\mathbb{F}$ for the cases $Char(\mathbb{F})=2$ and $Char(\mathbb{F})=3$ can be obtained by using Theorems \ref{Char2} and \ref{Char3}, respectively. The proof is similar to that of Theorem \ref{char0Endo}. Therefore, we give here the results without proof.

\begin{theor} Any non-trivial $2$-dimensional endo-commutative algebra over a field $\mathbb{F}$ $(Char(\mathbb{F})=2)$ is isomorphic to only one of the following, listed by their matrices of structure constants, such algebras:
\begin{itemize}
	\item $ A_{2,2}(\alpha_1,0,1)=\begin{pmatrix}
	\alpha_1&0&0&0 \\
	1& 1& 1+\alpha_1& 0
	\end{pmatrix}
	,$ where $ \alpha_1\in\mathbb{F}$
	
\item $ A_{3,2}(\alpha_1,0,\beta_2)=\begin{pmatrix}
\alpha_1&0&0&0 \\
0& \beta_2& 1+\alpha_1& 0
\end{pmatrix}
,$ where $ \alpha_1,\beta_2\in\mathbb{F}$

\item $ A_{4,2}(\alpha_1,\beta_1,0)=\begin{pmatrix}
\alpha_1&1&1&0 \\
\beta_1& 0& 1+\alpha_1& 1
\end{pmatrix}\simeq\begin{pmatrix}
\alpha_1&1&1&0 \\
\beta_1+a+a^2& 0& 1+\alpha_1& 1
\end{pmatrix}
,$\\ where $a, \alpha_1, \beta_1\in\mathbb{F}$

\item $ A_{5,2}(1,0)=\begin{pmatrix}
1&0&0&0 \\
1& 0& 1& 0
\end{pmatrix}
$

\item
$ A_{6,2}(\alpha_1,0)=\begin{pmatrix}
\alpha_1&0&0&0 \\
0& 1+\alpha_1& \alpha_1& 0
\end{pmatrix},
$ where $ \alpha_1\in\mathbb{F}$

\item $ A_{7,2}(0,\beta_1)=\begin{pmatrix}
0&1&1&0 \\
\beta_1& 1& 0& 1
\end{pmatrix}\simeq\begin{pmatrix}
0&1&1&0 \\
\beta_1+a+a^2& 1& 0& 1
\end{pmatrix}
,$ where $\beta_1, a\in\mathbb{F}$

\item $ A_{8,2}(\beta_1)=\begin{pmatrix}
0&1&1&1 \\
\beta_1& 0&0& 1
\end{pmatrix}\simeq\begin{pmatrix}
0&1&1&1 \\
\beta_1^{'}(a)& 0&0& 1
\end{pmatrix}
,$ where the polynomial $(\beta_1t^3+t+1)(\beta_1t^2+\beta_1t+1)$ has no root in $\mathbb{F},$
 $a\in\mathbb{F}$ and $ \beta_1^{'}(t)=\frac{(\beta_1^2t^3+\beta_1t+\beta_1)^2}{(\beta_1t^2+\beta_1t+1)^3}$

\item $ A_{9,2}(\beta_1)=\begin{pmatrix}
0&0&0&1\\
\beta_1&0&0&0
\end{pmatrix}\simeq\begin{pmatrix}
0&0&0&1\\
a^3\beta_1^{\pm1}&0&0&0
\end{pmatrix}
,$ where $\beta_1, a\in\mathbb{F},$ $a\neq 0$ and \\
the polynomial $ \beta_1+t^3$ has no root in $\mathbb{F}$

\item $ A_{10,2}(\beta_1)=\begin{pmatrix}
1&1&1&0\\
\beta_1&1&1&1
\end{pmatrix}\simeq\begin{pmatrix}
1&1&1&0\\
\beta_1+a+a^2&1&1&1
\end{pmatrix}
,$ where $ \beta_1, a\in\mathbb{F}$

\item $ A_{11,2}(\beta_1)=\begin{pmatrix}
0&1&1&0 \\
\beta_1&0&0&1\\
\end{pmatrix}\simeq\begin{pmatrix}
0&1&1&0 \\
b^2(\beta_1+a^2)&0&0&1\\
\end{pmatrix}
,$ where $a,b\in\mathbb{F}$ and $ b\neq 0$

\item $ A_{12,2}=\begin{pmatrix}
0&0&0&0 \\
1&0&0&0
\end{pmatrix}.
$
\end{itemize}

If $\mathbb{F}=\mathbb{Z}_2$ we obtain the list of representatives of endo-commutative algebras over $\mathbb{Z}_2$, as a particular case of the theorem, as follows.
\begin{corollary} \label{cor} Any non-trivial $2$-dimensional endo-commutative algebra over a field $\mathbb{F}=\mathbb{Z}_2,$ is isomorphic to only one of the following listed, by their matrices of structure constants, such algebras:

 $\bullet$ $ A_{2,2}(0,0,1)=\begin{pmatrix}
	0&0&0&0 \\
	1& 1& 1& 0
	\end{pmatrix}$ \ \
$\bullet$ $A_{2,2}(1,0,1)=\begin{pmatrix}
	1&0&0&0 \\
	1& 1& 0& 0
	\end{pmatrix} \ \
	$
$\bullet$ $ A_{3,2}(0,0,0)=\begin{pmatrix}
	0&0&0&0 \\
	0& 0& 1& 0
	\end{pmatrix}$

$\bullet$ $ A_{3,2}(0,0,1)=\begin{pmatrix}
	0&0&0&0 \\
	0& 1& 1& 0
	\end{pmatrix}$ \ \
$\bullet$
	$A_{3,2}(1,0,0)=\begin{pmatrix}
	1&0&0&0 \\
	0& 0& 0& 0
	\end{pmatrix}$ \ \
$\bullet$ $A_{3,2}(1,0,1)=\begin{pmatrix}
	1&0&0&0 \\
	0& 1&0& 0
	\end{pmatrix}$
		
 $\bullet$ $ A_{4,2}(0,0,0)=\begin{pmatrix}
	0&1&1&0 \\
	0& 0& 1& 1
	\end{pmatrix}$ \ \
$\bullet$ $	A_{4,2}(0,1,0)=\begin{pmatrix}
	0&1&1&0 \\
	1& 0& 1& 1
	\end{pmatrix}$ \ \
$\bullet$ $ A_{4,2}(1,0,0)=\begin{pmatrix}
	1&1&1&0 \\
	0& 0&0& 1
	\end{pmatrix}$

$\bullet$	$A_{4,2}(1,1,0)=\begin{pmatrix}
	1&1&1&0 \\
	1& 0&0& 1
	\end{pmatrix}$ \ \
		$\bullet$ $ A_{5,2}(1,0)=\begin{pmatrix}
	1&0&0&0 \\
	1&0&1& 0
	\end{pmatrix}
	$ \ \ \ \ \
 $\bullet$ $ A_{6,2}(0,0)=\begin{pmatrix}
	0&0&0&0 \\
	0& 1& 0& 0
	\end{pmatrix}$

	$\bullet$ $A_{6,2}(1,0)=\begin{pmatrix}
	1&0&0&0 \\
	0&0&1& 0
	\end{pmatrix}$ \ \
		\ \ \ $\bullet$ $ A_{7,2}(0,0)=\begin{pmatrix}
	0&1&1&0 \\
	0& 1& 0& 1
	\end{pmatrix}
	$ \ \ $\bullet$ $ A_{7,2}(0,1)=\begin{pmatrix}
	0&1&1&0 \\
	1& 1& 0& 1
	\end{pmatrix}
	$

	$\bullet$ $ A_{8,2}(1)=\begin{pmatrix}
	0&1&1&1 \\
	1& 0&0& 1
	\end{pmatrix}$ \ \
	\ \
\ \ \ \ $\bullet$ $ A_{10,2}(0)=\begin{pmatrix}
	1&1&1&0\\
	0&1&1&1
	\end{pmatrix}$ \ \
	\ 	$\bullet$ $ A_{10,2}(1)=\begin{pmatrix}
	1&1&1&0\\
	1&1&1&1
	\end{pmatrix}
	$
	
	$\bullet$ $A_{11,2}(0)=\begin{pmatrix}
	0&1&1&0 \\
	0&0&0&1\\
	\end{pmatrix}
	 $ \ \
	\quad \ \ $\bullet$ $A_{12,2}=\begin{pmatrix}
	0&0&0&0 \\
	1&0&0&0
	\end{pmatrix}.$

\end{corollary}

\end{theor}

\begin{theor} Any non-trivial two-dimensional endo-commutative algebra over a field $\mathbb{F}$ $(Char(\mathbb{F})=3)$ is isomorphic to only one of the following listed, by their matrices of structure constants, such algebras:
\begin{itemize}
\item
$A_{3,3}(\alpha_1,0,\beta_2)=\begin{pmatrix}
\alpha_1&0&0&0 \\
0& \beta_2& 1-\alpha_1& 0
\end{pmatrix},
\ \mbox{where} \ \alpha_1,\beta_2\in\mathbb{F}$
\item $A_{3,3}(2,\alpha_4, 1)=\begin{pmatrix}
2&0&0&\alpha_4 \\
0& 1& 2& 0
\end{pmatrix}\simeq\begin{pmatrix}
2&0&0&a^2\alpha_4 \\
0& 1& 2& 0
\end{pmatrix},
\ \mbox{where}\ a, \alpha_4\in\mathbb{F}$ and $a, \alpha_4\neq 0$
\item $A_{3,3}(2,\alpha_4, 2)=\begin{pmatrix}
2&0&0&\alpha_4 \\
0& 2& 2& 0
\end{pmatrix}\simeq\begin{pmatrix}
2&0&0&a^2\alpha_4 \\
0& 2& 2& 0
\end{pmatrix},
\ \mbox{where}\ a, \alpha_4\in\mathbb{F}$ and $a, \alpha_4\neq 0$
\item
$A_{5,3}(\alpha_1)=\begin{pmatrix}
\alpha_1&0&0&0 \\
1& 2\alpha_1-1& 1-\alpha_1& 0
\end{pmatrix},
\ \mbox{where}\ \alpha_1\in\mathbb{F}$
\item $ A_{7,3}(\alpha_1,0)=\begin{pmatrix}
\alpha_1&0&0&0 \\
0& 1-\alpha_1& -\alpha_1& 0
\end{pmatrix},
\ \mbox{where} \ \alpha_1\in\mathbb{F}$
\item $A_{7,3}(2,\alpha_4)=\begin{pmatrix}
2&0&0&\alpha_4 \\
0& 2& 1& 0
\end{pmatrix}\simeq\begin{pmatrix}
2&0&0&a^2\alpha_4 \\
0& 2& 1& 0
\end{pmatrix},$ where $a, \alpha_4\in\mathbb{F},$ and $a, \alpha_4\neq 0$
\item $A_{9,3}(\beta_1)=\begin{pmatrix}
0&1&1&1\\
\beta_1&0&0&-1
\end{pmatrix}\simeq\begin{pmatrix}
0&1&1&1\\
\beta'_1(a)&0&0&-1
\end{pmatrix},
$  where the polynomial\\ $(\beta_1-t^3)(\beta_1t^2+\beta_1t+1)(\beta_1^2t^3+\beta_1-2)$ has no root in $\mathbb{F},$
 $a\in\mathbb{F}$ and $ \beta_1^{'}(t)=\frac{(\beta_1^2t^3+\beta_1-2)^2}{(\beta_1t^2+\beta_1t+1)^3}$
\item $A_{10,3}(\beta_1)=\begin{pmatrix}
0&0&0&1\\
\beta_1&0&0&0
\end{pmatrix}\simeq\begin{pmatrix}
0&0&0&1\\
a^{3}\beta_1^{\pm1}&0&0&0
\end{pmatrix},
$  where polynomial $ \beta_1-t^3$ has no root, $a, \beta_1\in\mathbb{F},$ and $a, \beta_1\neq 0$

\item $ A_{11,3}(\beta_1)=\begin{pmatrix}
0&1&1&0 \\
\beta_1&0&0&-1\\
\end{pmatrix}\simeq\begin{pmatrix}
0&1&1&0 \\
a^{2}\beta_1&0&0&-1\\
\end{pmatrix},
\ \mbox{where}\  a, \beta_1\in\mathbb{F},$ and $a\neq 0,$
\item $ A_{12,3}=\begin{pmatrix}
1&0&0&0 \\
1&-1&-1&0
\end{pmatrix}
$
\item $ A_{13,3}=\begin{pmatrix}
0&0&0&0 \\
1&0&0&0
\end{pmatrix}.
$
\end{itemize}

\end{theor}

\section{Two-dimensional curled algebras} \label{CCA}
%

Let $\mathbb{A}$ be a two-dimensional algebra and  \[A=\left(\begin{array}{cccc} \alpha_1^{'} & \alpha_2^{'} & \alpha_3^{'} &\alpha_4^{'}\\ \beta_1^{'} & \beta_2^{'} & \beta_3^{'} &\beta_4^{'}\end{array}\right)\] be its MSC with respect to a basis $\mathbf{e}=\{\mathrm{e}_1, \mathrm{e}_2\}$. Then (\ref{C}) in terms of the elements of $A$ and the coordinate vector of $\mathrm{x}$ is written as follows
\begin{equation} \label{CA} \left\lbrace\begin{array}{rrl}
\alpha_1^{'}x_1^2+(\alpha_2^{'}+\alpha_3^{'})x_1x_2+\alpha_4^{'}x_2^2-\lambda(x_1,x_2)x_1&=&0\\
\beta_1^{'}x_1^2+(\beta_2^{'}+\beta_3^{'})x_1x_2+\beta_4^{'}x_2^2-\lambda(x_1,x_2)x_2&=&0.\\
\end{array}\right.
\end{equation}

We consider two options:
\begin{enumerate}
	\item Let $\mathrm{\lambda(x)}$ be not identically zero over $ \mathbb{F}^2\setminus\{(0,0)\}$. Since  (\ref{CA}) is an identity it must hold for any $x=(x_1,x_2)\in \mathbb{F}^2.$ Particularly, it must hold for $x=(0,x_2)$ and $x=(x_1,0)$ as well. Therefore, $x_1=0, x_2\neq 0$ ($x_1\neq 0, x_2=0$) implies
$\alpha_4^{'}=0,\ \lambda(0,x_2)=\beta_4^{'}x_2$  (respectively, $\beta_1^{'}=0,\ \lambda(x_1,0)=\alpha_1^{'}x_1$) and (\ref{CA}) becomes
$$ \left\lbrace\begin{array}{rrl}
\alpha_1^{'}x_1^2+(\alpha_2^{'}+\alpha_3^{'})x_1x_2-\lambda(x_1,x_2)x_1&=&0\\
(\beta_2^{'}+\beta_3^{'})x_1x_2+\beta_4^{'}x_2^2-\lambda(x_1,x_2)x_2&=&0.
\end{array}\right.
$$
Thus, if
$x_1\neq 0,\ x_2\neq 0$ then the system (\ref{CA}) can be rewritten as follows
$$ \left\lbrace\begin{array}{rrl}
\alpha_1^{'}x_1+(\alpha_2^{'}+\alpha_3^{'})x_2-\lambda(x_1,x_2)&=&0\\
(\beta_2^{'}+\beta_3^{'})x_1+\beta_4^{'}x_2-\lambda(x_1,x_2)&=&0.
\end{array}\right.
$$ This gives $ (\alpha_1^{'}-\beta_2^{'}-\beta_3^{'})x_1+(\alpha_2^{'}+\alpha_3^{'}-\beta_4^{'})x_2=0$ and, if  $\mathrm{Card}(\mathbb{F})>2$ (that is $\mathbb{F}\neq\mathbb{Z}_2$) then the later is noting else $\alpha_1^{'}=\beta_2^{'}+\beta_3^{'},\ \beta_4^{'}=\alpha_2^{'}+\alpha_3^{'}$. So, if $\mathbb{F}\neq\mathbb{Z}_2$ then  we have $\alpha_1^{'}x_1+\beta_4^{'}x_2=\lambda(x_1,x_2),$ at least one of $ \alpha_1^{'},\beta_4^{'}$ is not zero and
\begin{equation} \label{CR1A} A=\left(\begin{array}{cccc} \alpha_1^{'} & \alpha_2^{'} & \beta_4^{'}-\alpha_2^{'} &0\\ 0 & \beta_2^{'} & \alpha_1^{'}-\beta_2^{'} &\beta_4^{'}\end{array}\right).\end{equation}

If $\mathbb{F}=\mathbb{Z}_2$ then $\alpha_1^{'}x_1+(a-\beta_2^{'}-\beta_3^{'})x_2=\lambda(x_1,x_2),$ at least one of $ \alpha_1^{'}$, $a-\beta_2^{'}-\beta_3^{'}$ is not zero and one has
\begin{equation} \label{CR2A} A=\left(\begin{array}{cccc} \alpha_1^{'} & \alpha_2^{'} & a-\alpha_1^{'}-\alpha_2^{'} &0\\ 0 & \beta_2^{'} & \beta_3^{'}&a-\beta_2^{'}-\beta_3^{'}\end{array}\right),\end{equation}
 where $a=\lambda(1,1)$.

\item
Let $\mathrm{\lambda(x)}$ be identically zero over $ \mathbb{F}^2\setminus\{(0,0)\}$. In this case the system (\ref{CA}) is equivalent to $\alpha_1^{'}=\alpha_4^{'}=\beta_1^{'}=\beta_4^{'}=\alpha_2^{'}+\alpha_3^{'}=\beta_2^{'}+\beta_3^{'}=0$ and therefore, we have
\begin{equation} \label{CR3A} A=\left(\begin{array}{cccc} 0 & \alpha_2^{'} & -\alpha_2^{'} &0\\ 0 & \beta_2^{'} & -\beta_2^{'} &0\end{array}\right).\end{equation}
These type of algebras were said to be zeropotent (see \cite{CedilnikJerman, KSTT1, KSTT, TST2}).\\
The forms (\ref{CR2A}), (\ref{CR2A}) and (\ref{CR3A}) are the conditions for $A$ to be curled.
\end{enumerate}

Since now we know for a two-dimensional algebra given by MSC $A$ to be a curled algebra, we can treat the description of two-dimensional curled algebras by using Theorems 1, 2 and 3. For example, in the case of $\mathbb{F}$ with $Char(\mathbb{F})\neq 2, 3,$ a two-dimensional algebra given by
 \[ A=A_1(\mathrm{c})=\begin{pmatrix}
\alpha_1&\alpha_2&1+\alpha_2&\alpha_4 \\
\beta_1& -\alpha_1& 1-\alpha_1& -\alpha_2
\end{pmatrix},
\] where $ \mathrm{c}=(\alpha_1,\alpha_2,\alpha_4,\beta_1)\in\mathbb{F}^4,$ is a curled algebra if and only if $A$ has the following form
\[A=\left(\begin{array}{cccc} \alpha_1^{'} & \alpha_2^{'} & \beta_4^{'}-\alpha_2^{'} &0\\ 0 & \beta_2^{'} & \alpha_1^{'}-\beta_2^{'} &\beta_4^{'}\end{array}\right),\] that is
$$ \left\lbrace\begin{array}{rllrrrr}

\alpha_1^{'}&=&\alpha_1\ \ \

&\alpha_2^{'}&=&\alpha_2\\

\beta_4^{'}-\alpha_2^{'}&=&1+\alpha_2\ \ \

&0&=&\alpha_4\\

0&=&\beta_1\ \ \

&\beta_2^{'}&=&-\alpha_1\\

\alpha_1^{'}-\beta_2^{'}&=&1-\alpha_1\ \ \

&\beta_4^{'}&=&-\alpha_2

\end{array}\right.
$$
which happens if and only if
$\alpha_1=\frac{1}{3}, \alpha_2=-\frac{1}{3}, \alpha_4=0,\beta_1=0$. Therefore, among the algebras $A_1(\mathrm{c})$ only
$A_1(\frac{1}{3},-\frac{1}{3},0,0)$ is curled.

Going through $A_2(\mathrm{c})-A_{13}$ of Theorem \ref{Char0} in this manner  one comes to the following result.
\begin{theor} Any non-trivial $2$-dimensional curled algebra over a field $\mathbb{F}$ $(Char(\mathbb{F})\neq 2,3)$ is isomorphic to only one of the following, listed by their matrices of structure constants, such algebras$:$
	
	\begin{itemize}
		\item
		$ A_1(\frac{1}{3},-\frac{1}{3},0,0)=\begin{pmatrix}
		\frac{1}{3}&-\frac{1}{3}&\frac{2}{3}&0 \\
		0& -\frac{1}{3}& \frac{2}{3}& \frac{1}{3}
		\end{pmatrix}
		,$
		\item
		$ A_3(\alpha_1,0,2\alpha_1-1)=\begin{pmatrix}
		\alpha_1&0&0&0 \\
		0& 2\alpha_1-1& 1-\alpha_1& 0
		\end{pmatrix},
		$ where $\alpha_1\in\mathbb{F} $
		\item
		$ A_7(\frac{1}{3},0)=\begin{pmatrix}
		\frac{1}{3}&0&0&0 \\
		0& \frac{2}{3}& -\frac{1}{3}& 0
		\end{pmatrix}
		$.
	\end{itemize}
	
	Among the listed algebras only $ A_3(0,0,-1)$ is zeropotent.
\end{theor}

Similarly, in $Char(\mathbb{F})= 2,3$ cases the following results hold true.

\begin{theor} Any non-trivial $2$-dimensional curled algebra over a field $\mathbb{F}\neq\mathbb{Z}_2$ $(Char(\mathbb{F})=2)$ is isomorphic to only one of the following, listed by their matrices of structure constants, such algebras:
\begin{itemize}
\item $ A_{1,2}(1,1,0,0)=\begin{pmatrix}
1&1&0&0 \\
0& 1& 0& 1
\end{pmatrix}
$,

\item $ A_{3,2}(\alpha_1,0,1)=\begin{pmatrix}
\alpha_1&0&0&0 \\
0& 1& 1+\alpha_1& 0
\end{pmatrix}
,$ where $\alpha_1\in\mathbb{F}$

\item
$ A_{6,2}(1,0)=\begin{pmatrix}
1&0&0&0 \\
0& 0& 1& 0
\end{pmatrix}
$

\end{itemize}

and among them only $ A_{3,2}(0,0,1)$ is zeropotent.

\end{theor}

		\begin{theor} \label{} Any non-trivial $2$-dimensional curled algebra over a field $\mathbb{F}=\mathbb{Z}_2,$ is isomorphic to only one of the following listed, by their matrices of structure constants, such algebras:
			\begin{itemize}
				\item
				$ A_{1,2}(1,1,0,0)=\begin{pmatrix}
				1&1&0&0 \\
				0& 1& 0&1
				\end{pmatrix},
				$

				\item $ A_{3,2}(\alpha_1,0,1)=\begin{pmatrix}
				\alpha_1&0&0&0 \\
				0& 1& 1+\alpha_1& 0
				\end{pmatrix},
				$ where $ \alpha_1\in\mathbb{Z}_2$

				\item $ A_{4,2}(\alpha_1,0,0)=\begin{pmatrix}
				\alpha_1&1&1&0 \\
				0& 0& 1+\alpha_1& 1
				\end{pmatrix},
				$ where $ \alpha_1\in\mathbb{Z}_2$

				\item $A_{6,2}(1,0)=\begin{pmatrix}
				1&0&0&0 \\
				0& 0& 1& 0
				\end{pmatrix},
				$
				\item $A_{7,2}(0,0)=\begin{pmatrix}
				0&1&1&0 \\
				0& 1& 0& 1
				\end{pmatrix},
				$

				\item $A_{10,2}(0)=\begin{pmatrix}
				1&1&1&0\\
				0&1&1&1
				\end{pmatrix}$.
				
			Among these algebras only $A_{3,2}(0,0,1)$ is zeropotent. 	
					
			\end{itemize}
		\end{theor}

\begin{theor} Any non-trivial $2$-dimensional curled algebra over a field $\mathbb{F}$ $(Char(\mathbb{F})=3)$ is isomorphic to only one of the following listed, by their matrices of structure constants, such algebras:
\begin{itemize}
\item
$A_{3,3}(\alpha_1,0,2+2\alpha_1)=\begin{pmatrix}
\alpha_1&0&0&0 \\
0& 2+2\alpha_1& 1+2\alpha_1& 0
\end{pmatrix},
\ \mbox{where} \  \ \alpha_1\in\mathbb{F},$
\item
$A_{4,3}(0,2)=\begin{pmatrix}
0&1&1&0 \\
0& 2& 1& 2
\end{pmatrix},
$
\item $ A_{11,3}(0)=\begin{pmatrix}
0&1&1&0 \\
0&0&0&2\\
\end{pmatrix},
$
\end{itemize}

Among the listed algebras only $ A_{3,3}(0,0,2)$ is zeropotent.

\end{theor}

\section{Two-dimensional endo-commutative curled algebras}
In this section we make use the results of the last two sections to get a classification of two-dimensional endo-commutative curled algebras up to isomorphism.

\begin{theor} Any non-trivial two-dimensional endo-commutative curled  algebra over a field $\mathbb{F}$ $(Char(\mathbb{F})\neq 2,3)$ is isomorphic to only one of the following, listed by their matrices of structure constants, such algebras:

\begin{itemize}
\item
$ A_3(\alpha_1,0,2\alpha_1-1)=\begin{pmatrix}
\alpha_1&0&0&0 \\
0& 2\alpha_1-1& 1-\alpha_1& 0
\end{pmatrix},
$ where $ \alpha_1\in\mathbb{F} ,$
\item
 $ A_7(\frac{1}{3},0)=\begin{pmatrix}
\frac{1}{3}&0&0&0 \\
0& \frac{2}{3}& -\frac{1}{3}& 0
\end{pmatrix}
$
\end{itemize}

and among them only $ A_3(0,0,-1)$ is zeropotent.
\end{theor}

\begin{theor} Any non-trivial $2$-dimensional Endo-commutative curled  algebra over a field $\mathbb{F}\neq \mathbb{Z}_2$ $(Char(\mathbb{F})=2)$ is isomorphic to only one of the following, listed by their matrices of structure constants, such algebras:
\begin{itemize}
\item $ A_{3,2}(\alpha_1,0,1)=\begin{pmatrix}
\alpha_1&0&0&0 \\
0& 1& 1+\alpha_1& 0
\end{pmatrix}
,$ where $\alpha_1\in\mathbb{F} ,$
\item
$ A_{6,2}(1,0)=\begin{pmatrix}
1&0&0&0 \\
0& 0& 1& 0
\end{pmatrix}
$.
\end{itemize}

Among the algebras listed above only $ A_{3,2}(0,0,1)$ is a zeropotent algebra.
\end{theor}.

\begin{theor} \label{13} Any non-trivial two-dimensional endo-commutative curled  algebra over a field $\mathbb{F}=\mathbb{Z}_2$ is isomorphic to only one of the following listed, by their matrices of structure constants, such algebras:

			$\bullet$ $ A_{3,2}(0,0,1)=\begin{pmatrix}
			0&0&0&0 \\
			0& 1& 1& 0
			\end{pmatrix}$ \ \
			$\bullet$
			$A_{3,2}(1,0,1)=\begin{pmatrix}
			1&0&0&0 \\
			0& 1& 0& 0
			\end{pmatrix}$ \ \
			$\bullet$ $A_{4,2}(0,1,0)=\begin{pmatrix}
			0&1&1&0 \\
			0& 0&1& 1
			\end{pmatrix}$

			$\bullet$ $ A_{4,2}(1,1,0)=\begin{pmatrix}
			1&1&1&0 \\
			0& 0& 0& 1
			\end{pmatrix}$ \ \
				$\bullet$ $ A_{6,2}(1,0)=\begin{pmatrix}
				1&0&0&0 \\
				0& 0& 1& 0
				\end{pmatrix}
				$

						$\bullet$ $ A_{7,2}(0,0)=\begin{pmatrix}
			0&1&1&0 \\
			0& 1& 0& 1
			\end{pmatrix}
			$ \ \ \quad
			$\bullet$ $ A_{10,2}(0)=\begin{pmatrix}
			1&1&1&0\\
			0&1&1&1
			\end{pmatrix}.$
		
Here, only $ A_{3,2}(0,0,1)$ is zeropotent.			
\end{theor}
Note that the endo-commutative curled algebras given in \cite{TST} can be found in Theorem \ref{13} as $ A_{3,2}(0,0,1)\simeq C_4, \ A_{3,2}(1,0,1)\simeq C_3,\ A_{4,2}(0,1,0)\simeq C_{12}, \ A_{4,2}(1,1,0)\simeq C_7,\ A_{6,2}(1,0)\simeq C_2, \ A_{7,2}(0,0)\simeq C_{13}$ and $A_{10,2}(0)\simeq C_1.$
\begin{theor} Any non-trivial two-dimensional endo-commutative curled  algebra over a field $\mathbb{F}$ $(Char(\mathbb{F})=3)$ is isomorphic to only one of the following listed, by their matrices of structure constants, such algebras:
\begin{itemize}
\item
$A_{3,3}(\alpha_1,0,2\alpha_1-1)=\begin{pmatrix}
\alpha_1&0&0&0 \\
0& 2\alpha_1-1& 1-\alpha_1& 0
\end{pmatrix},
$
where $\alpha_1\in\mathbb{F},$
\item
$A_{11,3}(0)=\begin{pmatrix}
0&1&1&0 \\
0& 0&0& 2
\end{pmatrix},
$


\end{itemize}

Among the algebras listed above only $A_{3,3}(0,0,-1)$ is zeropotent.

\end{theor}

\section{Acknowledgement}
We would like to thank Professor K.Shirayanagi for pointing out the authors' attention to an algebra that was missed in earlier version of Corollary \ref{cor}.

\end{document}